\documentclass[11pt,reqno]{amsart}
\usepackage[a4paper,margin=28mm]{geometry}
\usepackage[T1]{fontenc}
\usepackage{lmodern}
\usepackage{amsmath,amssymb,amsthm,mathtools,mathrsfs}
\usepackage{enumitem}
\usepackage{hyperref}
\hypersetup{
  colorlinks=true,
  linkcolor=blue!55!black,
  citecolor=blue!55!black,
  urlcolor=blue!55!black,
  pdftitle={Ramification of torsion packets on curves with good-reduction Jacobians},
  pdfauthor={Arata Minamide}
}
\usepackage{xcolor}
\theoremstyle{plain}
\newtheorem{theorem}{Theorem}[section]
\newtheorem{proposition}[theorem]{Proposition}
\newtheorem{lemma}[theorem]{Lemma}

\newtheorem*{conjecture}{Conjecture}
\newtheorem*{maintheorem}{Main Theorem}
\newtheorem*{introcorollary}{Corollary}
\theoremstyle{definition}

\newtheorem*{setting}{Setting}
\newtheorem*{fact}{Fact}
\theoremstyle{remark}

\newcommand{\bbF}{\mathbb F}
\newcommand{\bP}{\mathbb P}

\newcommand{\cH}{\mathcal H}

\newcommand{\cJ}{\mathcal J}
\newcommand{\cD}{\mathcal D}

\newcommand{\Kbar}{\overline{K}}

\newcommand{\divisor}{\operatorname{div}}

\newcommand{\Gal}{\operatorname{Gal}}

\newcommand{\Fp}{\bbF_p}
\title[Ramification of torsion packets on curves with good-reduction Jacobians]{Ramification of torsion packets on curves with good-reduction Jacobians}
\author{Arata Minamide}
\address{ZEN Mathematics Center, ZEN University, Zushi 249-0007, Japan; Research Institute for Mathematical Sciences, Kyoto University, Kyoto 606-8502, JAPAN}
\email{arata\_minamide@zen.ac.jp; minamide@kurims.kyoto-u.ac.jp}
\date{November 2026}
\begin{document}
\begin{abstract}
Let $p\ge5$ be a prime number, $k$ an algebraically closed field of characteristic $p$. Write $R=W(k)$ for the ring of Witt vectors over $k$, and $K$ for the field of fractions of $R$. Let $X$ be a smooth, proper, geometrically connected curve over $K$ of genus greater than one. Write $J$ for the Jacobian variety of $X$. Fix an algebraic closure $\Kbar$ of $K$. In the present paper, we prove that if $J$ has good reduction over $R$, then every point $x\in X(\Kbar)$ such that $[\sigma x-x]$ is torsion for each $\sigma\in \operatorname{Gal}(\Kbar/K)$ is $K$-rational.
\end{abstract}
\maketitle
\tableofcontents

\setcounter{section}{-1}
\section{Introduction}\label{sec:intro}
In \cite{Coleman}, R. Coleman conjectured that Galois-stable torsion packets on curves over number fields are unramified at certain primes of good reduction [cf. \cite{Coleman}, Conjecture B; \cite{ColemanKaskelRibet}, Errata].

In the present paper, we consider the following local assertion.

\begin{setting}
Let $p\ge5$ be a prime number, $k$ an algebraically closed field of characteristic $p$, $R=W(k)$ the ring of Witt vectors over $k$, and $K$ the field of fractions of $R$. Fix an algebraic closure $\Kbar$ of $K$, and write $G_K=\Gal(\Kbar/K)$.
Let $X$ be a smooth, proper, geometrically connected curve over $K$ of genus $g\ge2$. 
Write $J$ for the Jacobian variety of $X$.
\end{setting}

\begin{conjecture}
Let $x\in X(\Kbar)$.
Suppose that 
\begin{itemize}
\item $J$ has good reduction over $R$, i.e., there exists an abelian scheme $\mathcal{J}$ over $R$ equipped with an isomorphism between its generic fiber $\mathcal{J}_K$ and $J$, and
\item for each $\sigma\in G_K$, $[\sigma x-x]\in J(\Kbar)$ is torsion.
\end{itemize}
Then $x\in X(K)$.
\end{conjecture}

\begin{fact}
We consider an additional assumption that
\begin{itemize}
\item[($\ast$)] there exists $x_0 \in X(K)$ such that $[x-x_0]$ is torsion.
\end{itemize}
The above Conjecture holds in the following cases:
\begin{enumerate}[label=\textup{(\roman*)},leftmargin=*,itemsep=3pt,topsep=3pt]
\item \textbf{Coleman (1987).}
$k=\overline{\mathbb F}_p$, $X$ has good reduction over $R$, and either $p>2g$ or the special fiber $\cJ_k$ is ordinary or superspecial [cf. \cite{Coleman}, Corollary 20.2].

\item \textbf{Tamagawa (2001).}
($\ast$) is satisfied, and $\cJ_k$ is ordinary [cf. \cite{Tamagawa}, Proposition 2.1 and Theorem 3.6, (ii)].

\item \textbf{Hoshi (2017).}
($\ast$) is satisfied, $g\ge p$, and either $X_{\overline K}=X\times_K\Kbar$ has gonality $2$ or $>p$, or all its Weierstrass points are $K$-rational [cf. \cite{Hoshi2017}, Theorem C].

\item \textbf{Hoshi (2022).}
($\ast$) is satisfied, and $\cJ_k$ is superspecial
[cf. \cite{Hoshi2022}, Theorem E].

\item \textbf{Takao (2023).}
($\ast$) is satisfied, $g=2$, and $\cJ_k$ is supersingular
[cf. \cite{Takao2023}, Theorem 3].
\end{enumerate}
\end{fact}

\begin{maintheorem}
The above Conjecture holds.
\end{maintheorem}

\begin{introcorollary}
Conjecture B of \cite{Coleman} [cf. \cite{ColemanKaskelRibet}, Errata] holds.
\end{introcorollary}

\begin{proof}
In the present proof, we use the notation of
\cite{Coleman}, Conjecture B, with the correction given in
\cite{ColemanKaskelRibet}, Errata.
Write $p$ for the residue characteristic of $\mathfrak p$, and 
\[
L \ \stackrel{\mathrm{def}}{=} \
\widehat{K_{\mathfrak p}^{\mathrm{ur}}} \ \ (= \ \widehat{\mathbb Q_p^{\mathrm{ur}}} \ \simeq \ W(\overline{\mathbb F}_p)[1/p])
\]
for the completion of the maximal unramified extension of $K_{\mathfrak p}$ [cf. the assumption that $K_{\mathfrak p}/\mathbb Q_p$ is unramified].
Moreover, since $C$ has good reduction at $\mathfrak p$,
the Jacobian of $C_L$ has good reduction over the ring
of integers of $L$.
Let $x\in T$.
Since $T$ is a Galois-stable torsion packet, it holds that $[\sigma x-x]$ is torsion for each $\sigma\in\Gal(\overline L/L)$.
Thus, our Main Theorem implies that $x\in C(L)$.
Since $x\in T$ was arbitrary, the inertia group at $\mathfrak p$ acts trivially on $T$.
Therefore, we conclude that $K(T)/K$ is unramified above $\mathfrak p$.
\end{proof}

\medskip
\noindent\textbf{Outline of the proof.}
Write $S=G_Kx$ and $m=|S|$.
Our goal is to prove that 
\[
m \ = \ 1.
\]
Suppose that $m>1$. Then we obtain two bounds for $m$ as follows.

The first bound is arithmetic. We consider the $\mathbb{Z}[G_K]$-submodule
\[
D \ \stackrel{\mathrm{def}}{=} \ \big\langle[y-x] \mid y \in S\big\rangle_{\mathbb{Z}}  \ \subseteq \ J(\Kbar)
\]
[which has essentially been considered in \cite{Hoshi2022}, Definition 2.3, (ii)].
If ($\ast$) is satisfied, then it follows from \cite{Hoshi2017}, Theorem 3.4, (ii), that $pD=\{0\}$. In our general setting, we prove this equality
by an argument similar to the original proof of \cite{Hoshi2017}, Theorem 3.4, (ii), but using an auxiliary finite flat group scheme.
In light of this equality, $D$ may naturally be regarded as an $\Fp[G_K]$-submodule of $J[p](\Kbar)$.
Then the schematic closure $\cD$ of $D$ in $\cJ[p]$ is a finite flat group scheme over $R$ with connected special fiber.
By applying Raynaud's character theory [cf. \cite{Raynaud}, Th\'eor\`eme 3.4.3]
to $\cD$, together with Guralnick's bound on cyclic quotients of transitive groups [cf. \cite{Guralnick2000}, Theorem 1] and ``cyclotomic displacement rigidity'', we obtain
\begin{equation}\label{eq:intro-lower}
m \ > \ p(p-1).
\end{equation}

The second bound is geometric. We consider the dimension
\[
h^0(px) \ \stackrel{\mathrm{def}}{=} \ \dim_{\Kbar} H^0(X_{\Kbar}, \mathcal{O}_{X_{\Kbar}}(px)) \ \ge \ 2
\]
[cf. the equality $pD=\{0\}$].
We first show that $h^0(px)\ge3$.
Indeed, if $h^0(px)=2$, then an augmentation-module inequality arising from the self-duality of the permutation module $\Fp[S]$ and Raynaud's character theory, together with a certain Kummer-theoretic restriction argument, forces $m=3$ and $\dim_{\Fp}D=1$.
Then a simpler variant of the character-theoretic argument used to obtain \eqref{eq:intro-lower} yields a contradiction.
Under the assumption that $h^0(px)\ge3$, the Pl\"ucker formula [cf. \cite{Laksov1984},
Theorems 9 and 15] and the classical bounds of Clifford, Castelnuovo, and Riemann--Roch
[cf. \cite{ACGH1985}, Chapter III] give
\begin{equation}\label{eq:intro-upper}
m \ \le \ (p-1)^2.
\end{equation}
Since $(p-1)^2<p(p-1)$, the bounds \eqref{eq:intro-lower} and \eqref{eq:intro-upper}
are incompatible. Thus, $m=1$, as required.

\medskip
\noindent\textbf{Organization of the paper.}
In \S 1, we establish rigidity results for finite flat group schemes over $R$.
In \S 2, we introduce $D$, prove that $pD=\{0\}$, and show that its schematic closure has connected special fiber.
In \S 3, we prove an orbit-augmentation inequality.
In \S 4, we combine the inequality of \S 3 with a Kummer-theoretic argument to show that $h^0(px)\ge3$ if $m>1$.
In \S 5, we exclude small nontrivial Galois orbits.
In \S 6, we establish a purely geometric bound for prime-degree linear systems on curves.
Finally, in \S 7, we prove our Main Theorem.

\section{Finite flat group schemes over an absolutely unramified base}\label{sec:finiteflat}
In the present section, we establish several rigidity results for finite flat group schemes.
Our main ingredients are Raynaud's classification of characters and Cartier duality.

In the present section, let
\begin{itemize}
\item $p\ge5$ be a prime number, and
\item $k$ an algebraically closed field of characteristic $p$.
\end{itemize}
Write
\begin{itemize}
\item $R$ for the ring of Witt vectors over $k$, 
\item $K$ for the field of fractions of $R$,
\item $G_K \stackrel{\mathrm{def}}{=} \operatorname{Gal}(\overline K/K)$ for the absolute Galois group of $K$, where $\overline K$ is a fixed algebraic closure of $K$, and
\item $P_K \subseteq G_K$ for the wild inertia subgroup.
\end{itemize}
For each positive integer $f$, write
\[
\omega_f: G_K \ \rightarrow \ \mathbb{F}^{\times}_{p^f} \ \subseteq \ \overline{\mathbb{F}}^{\times}_p \ \ (\subseteq \ k^{\times})
\]
for the composite of 
\begin{itemize}
\item[-] the natural surjection $G_K \twoheadrightarrow G_K/P_K$,
\item[-] the natural isomorphism $G_K/P_K \stackrel{\sim}{\rightarrow} \varprojlim_n \mathbb F_{p^n}^{\times}$, where, for any $n_1\mid n_2$, the transition homomorphism $\mathbb F_{p^{n_2}}^{\times} \twoheadrightarrow \mathbb F_{p^{n_1}}^{\times}$ is given by sending
\[
a \ \mapsto \ a^{(p^{n_2}-1)/(p^{n_1}-1)}
\]
[cf. \cite{Raynaud}, \S 3.1; also \cite{Serre1972}, \S 1.3, Proposition 2], and 
\item[-] the natural projection $\varprojlim_n \mathbb F_{p^n}^{\times} \twoheadrightarrow \mathbb F_{p^f}^{\times}$.
\end{itemize}
We note that $\omega_1 = \omega^{(p^f-1)/(p-1)}_f$ coincides with 
\begin{itemize}
\item the mod-$p$ cyclotomic character $\chi^{\mathrm{cyc}}_p:G_K\rightarrow\mathbb F_p^\times$
\end{itemize}
[cf. \cite{Serre1972}, \S 1.8, Proposition 8 and its Corollary].

Let 
\begin{itemize}
\item $\mathcal{G}$ be a finite flat group scheme over $R$.
\end{itemize}
For each $R$-scheme $\mathcal Z$, write
\begin{itemize}
\item $\mathcal Z_K$ for the generic fiber of $\mathcal Z$,
\item $\mathcal Z_{\overline K}$ for the geometric generic fiber of $\mathcal Z$, and
\item $\mathcal Z_k$ for the special fiber of $\mathcal Z$.
\end{itemize}

\begin{lemma}[Finite flat schematic closure]\label{lem:ff-closure}
Let $\mathcal{H}_K \subseteq \mathcal{G}_K$ be a closed subgroup scheme. 
Then the schematic closure of $\mathcal{H}_K$ in $\mathcal{G}$ is a finite flat closed subgroup scheme of $\mathcal{G}$ over $R$.
\end{lemma}

\begin{proof}
This follows immediately from \cite{Raynaud}, \S 2.1.
\end{proof}

\begin{lemma}[Raynaud's character classification]\label{lem:raynaud}
Suppose that 
\begin{itemize}
\item $\mathcal{G}$ is commutative, and
\item $\mathcal{G}$ is annihilated by $p$.
\end{itemize}
Write $V \stackrel{\mathrm{def}}{=} \mathcal G_K(\Kbar)$. 
Let
\[
\{0\} \ = \ V_0 \ \subsetneq \ V_1 \ \subsetneq \ \cdots \ \subsetneq \ V_r \ = \ V\otimes_{\mathbb F_p}\overline{\mathbb F}_p
\]
be a Jordan--Hölder filtration of $V\otimes_{\mathbb F_p}\overline{\mathbb F}_p$
as an $\overline{\mathbb F}_p[G_K]$-module.
For each $i\in\{1,\ldots,r\}$, write $Q_i \stackrel{\mathrm{def}}{=} V_i/V_{i-1}$.
Then the following hold:
\begin{itemize}
\item[(i)] Each $Q_i$ is one-dimensional. 
In particular, we have $\dim_{\mathbb{F}_p}V=r$.
Write
\[
\chi_i:G_K \ \rightarrow \ \overline{\mathbb F}_p^\times
\]
for the associated character.

\item[(ii)] Suppose that $\mathcal G_k$ is connected.
Then each $\chi_i$ is nontrivial.
Moreover, there exists a positive integer $f$ such that, for each $i\in\{1,\ldots,r\}$,
\[
\chi_i \ = \ \omega_f^{E_i}, \qquad E_i \ = \ \sum_{j=0}^{f-1}\varepsilon_{i,j}p^j
\]
--- where $\varepsilon_{i,j} \in\{0,1\}$.

\item[(iii)] Suppose that $\mathcal G_k$ is connected. 
Then, for any [not necessarily distinct] $s, t\in\{1,\ldots,r\}$, it holds that
\[
\chi_s \ \ne \ \chi_t^{-1}.
\]
\end{itemize}
\end{lemma}

\begin{proof}
First, we consider assertion (i).
Write 
\[
\chi_i:G_K \ \rightarrow \ \operatorname{Aut}_{\overline{\mathbb F}_p}(Q_i)
\]
for the representation associated to the $\overline{\mathbb F}_p[G_K]$-module $Q_i$.
Since $P_K$ is pro-$p$ and $\chi_i$ has finite image, it follows from \cite{Serre}, Proposition 26, that $Q_i^{P_K}\ne 0$.
Moreover, since $Q_i^{P_K}$ is an $\overline{\mathbb F}_p[G_K]$-submodule of $Q_i$ [cf. the fact that $P_K$ is normal in $G_K$], it follows from the irreducibility of $Q_i$ that $Q_i^{P_K}=Q_i$, hence that $\chi_i(P_K)=\{1\}$.
Thus, $\chi_i$ factors through a representation 
\[
\overline{\chi}_i: G_K/P_K \ \rightarrow \ \operatorname{Aut}_{\overline{\mathbb F}_p}(Q_i).
\]
Then since $G_K/P_K$ is procyclic [cf. our assumption that $k$ is algebraically closed] and $\overline{\chi}_i$ has finite image, it follows from the fact that $\overline{\mathbb F}_p$ is algebraically closed that $Q_i$ contains a one-dimensional $\overline{\mathbb F}_p[G_K]$-submodule.
Therefore, we conclude from the irreducibility of $Q_i$ that $Q_i$ is one-dimensional.
This completes the proof of assertion (i).

Next, we consider assertion (ii).
Since $\{V_i\}_{0\le i \le r}$ is a finite set of finite-dimensional
$\overline{\mathbb F}_p$-subspaces of
$V\otimes_{\mathbb F_p}\overline{\mathbb F}_p$, there exist a finite subfield $F=\mathbb F_{p^f}\subseteq\overline{\mathbb F}_p$ and a finite set $\{U_i\}_{0\le i \le r}$ of $F$-subspaces of $V\otimes_{\mathbb F_p}F$ such that
\[
U_i\otimes_F\overline{\mathbb F}_p \ = \ V_i.
\]
In particular, for each $i\in\{0,\ldots,r\}$, we have 
\[
U_i \ =\ V_i \cap (V\otimes_{\mathbb F_p}F).
\]
In light of this equality, we obtain a filtration
\[
\{0\} \ = \ U_0 \ \subsetneq \ U_1 \ \subsetneq \ \cdots \ \subsetneq \ U_r \ = \ V\otimes_{\mathbb F_p}F
\]
by $F[G_K]$-submodules.
Here, we note that each $U_i/U_{i-1}$ is one-dimensional over $F$, and the associated character is $\chi_i$.

Observe that since $[p]_{\mathcal G}=0$, $\mathcal G$ may naturally be regarded as an $\mathbb F_p$-module object in the category of commutative $R$-group schemes.
Write
\[
\mathcal{A} \ \stackrel{\mathrm{def}}{=} \ F\otimes_{\mathbb F_p}\mathcal G
\]
for the commutative $R$-group scheme equipped with
a natural $F$-action obtained from $\mathcal{G}$ by the Serre tensor construction
[cf. \cite{ConradGZ}, Theorem 7.2].
Since $F$ is a free $\mathbb F_p$-module of rank $f$,
$\mathcal A$ is isomorphic to 
\[
\underbrace{\mathcal G\times_R\cdots\times_R\mathcal G}_{f\text{ factors}}
\]
as an $R$-group scheme [cf. the proof of \cite{ConradGZ}, Theorem 7.2].
In particular, $\mathcal A$ is finite flat over $R$, and $\mathcal A_k$ is connected.
Moreover, by the definition of the Serre tensor construction,
we have a natural isomorphism
\[
\mathcal A_K(\overline K) \  \stackrel{\sim}{\rightarrow} \ V\otimes_{\mathbb F_p}F
\]
of $F[G_K]$-modules.

For each $i\in\{0,\ldots,r\}$, write $\mathcal U_i$ for the schematic closure in $\mathcal{A}$ of the closed subgroup scheme of $\mathcal{A}_K$ corresponding to $U_i$.
By Lemma \ref{lem:ff-closure}, each \(\mathcal U_i\) is a closed subgroup scheme of $\mathcal A$ which is finite flat over $R$.
Since $U_i\subsetneq U_{i+1}$, it follows from the definition
of schematic closure that
$\mathcal U_i\subseteq\mathcal U_{i+1}$.
Moreover, since each $U_i$ is an $F$-subspace of $\mathcal A_K(\overline K)$, the natural $F$-action on $\mathcal A$ preserves $\mathcal U_i$.
Thus, we obtain a filtration
\[
e_{\mathcal A}(\operatorname{Spec}R) \ = \ \mathcal{U}_0 \ \subsetneq \ \mathcal{U}_1 \ \subsetneq \ \cdots \ \subsetneq \ \mathcal{U}_r \ = \ \mathcal{A}
\]
by finite flat closed $R$-subgroup schemes equipped with compatible $F$-actions.
Then, for each $i\in\{1,\ldots,r\}$, write
\[
\mathcal H_i \ \stackrel{\mathrm{def}}{=} \
\mathcal U_i/\mathcal U_{i-1}
\]
for the fppf quotient.
Here, we note that since $\mathcal U_{i-1}$ is a finite flat closed subgroup scheme of $\mathcal U_i$, $\mathcal H_i$ is represented by a finite flat commutative $R$-group scheme [cf. \cite{Raynaud}, \S 2.1].
Moreover, the $F$-action on $\mathcal U_i$ induces an $F$-action on $\mathcal H_i$.
Since the formation of the fppf quotient commutes with base change
[cf. \cite{SGA3}, Exposé VI$_{\mathrm B}$, Proposition 9.2 (v)], it follows from the fact that $H^1_{\mathrm{fppf}}(\overline{K}, (\mathcal{U}_{i-1})_{\overline{K}})=\{0\}$ that we have a natural isomorphism
\[
(\mathcal H_{i})_K(\overline K) \ \stackrel{\sim}{\rightarrow} \ U_i/U_{i-1}
\]
of $F[G_K]$-modules.
In particular, each $(\mathcal H_{i})_K(\overline K)$ is one-dimensional over $F$, and the associated character is $\chi_i$.

With these preparations, let us examine the characters $\chi_i$.
We begin by observing that since $\mathcal A_k$ is the spectrum of an Artinian local ring, the closed subscheme $(\mathcal U_i)_k \subseteq \mathcal A_k$ is connected.
In particular, since the quotient morphism $\mathcal U_i\rightarrow\mathcal H_i$
is faithfully flat, hence surjective, it follows that $(\mathcal H_i)_k$ is connected.
Now suppose that $\chi_i$ is trivial.
Then
\[
(\mathcal H_{i})_K \ \simeq \ (\mathbb Z/p\mathbb Z)_K^f.
\]
In light of our assumption that the absolute ramification index of $R$ is $1$ ($<p-1$), it follows from \cite{Raynaud}, Théorème 3.3.3, that
\[
\mathcal H_i \ \simeq \ (\mathbb Z/p\mathbb Z)_R^f.
\]
This contradicts the connectedness of $(\mathcal H_i)_k$.
Thus, $\chi_i$ is nontrivial.
Moreover, since $\chi_i$ is tamely ramified [as shown in the proof of assertion (i)], it follows from \cite{Raynaud}, Théorème 3.4.3, that
\[
\chi_i \ = \ \omega_f^{E_i}, \qquad E_i \ = \ \sum_{j=0}^{f-1}\varepsilon_{i,j}p^j,
\qquad \varepsilon_{i,j} \ \in \ \{0, 1\}.
\]
This completes the proof of assertion (ii).

Finally, we consider assertion (iii).
By assertion (ii), for each $i\in\{1,\ldots,r\}$, we have
\[
0 \ < \ E_i \ \le \ \sum_{j=0}^{f-1}p^j \ = \ \frac{p^f-1}{p-1},
\]
where the first inequality follows from the nontriviality of $\chi_i$.
Thus, for any [not necessarily distinct] $s, t\in\{1,\ldots,r\}$, we obtain that
\[
0 \ < \ E_s+E_t \ \le \ \frac{2(p^f-1)}{p-1} \ < \ p^f-1,
\]
where the last inequality follows from our assumption that $p\ge 5$.
Therefore, since $\omega_f$ has order $p^f-1$, we conclude that
\[
\chi_s\chi_t \ = \ \omega_f^{E_s+E_t} \ \neq \ 1,
\]
hence that $\chi_s\neq\chi_t^{-1}$.
This completes the proof of assertion (iii).
\end{proof}

\begin{lemma}[Low-order character rigidity]\label{lem:digit}
We continue to use the notation of Lemma \ref{lem:raynaud}. 
Let $i\in\{1,\ldots,r\}$.
Suppose that 
\begin{itemize}
\item $\mathcal G_k$ is connected, and
\item $\chi_i$ has order $d_i \le p(p-1)$.
\end{itemize}
Then $\chi_i$ is the mod-$p$ cyclotomic character.
\end{lemma}

\begin{proof}
Write
\[
N \ \stackrel{\mathrm{def}}{=} \ p^f-1, \qquad E_{\mathrm{cyc}} \ \stackrel{\mathrm{def}}{=} \ \frac{N}{p-1}.
\]
By Lemma \ref{lem:raynaud}, (ii), we have
\[
0 \ < \ E_i \ \le \ E_{\mathrm{cyc}}.
\]
Now suppose that $E_i\neq E_{\mathrm{cyc}}$.
Then the base-$p$ digits of $E_i$ contain both $0$ and $1$.
Since $p^f\equiv 1\pmod N$, multiplication by $p$ modulo $N$ cyclically permutes these digits.
Thus, after a suitable cyclic permutation, we obtain an exponent $E'$ satisfying
\[
p\mid E',
\qquad
0 \ < \ E' \ < \ E_{\mathrm{cyc}}.
\]
Moreover, since $p$ is relatively prime to $N$, the character $\omega_f^{E'}$ has the same order $d_i$ as $\chi_i$.
In particular, there exists a positive integer $\ell$ such that
\[
d_iE' \ = \ \ell N.
\]
Then our assumption that $d_i\le p(p-1)$ implies that
\[
0 \ < \ \ell \ = \ \frac{d_iE'}{N} \ < \ \frac{d_iE_{\mathrm{cyc}}}{N} \ = \ \frac{d_i}{p-1} \ \le \ p.
\]
On the other hand, since $N \equiv -1 \pmod p$, we obtain that
\[
0 \ \equiv \ d_iE' \ = \ \ell N \ \equiv \ -\ell \ \pmod p,
\]
a contradiction. Thus, we conclude that $E_i=E_{\mathrm{cyc}}$. 

Therefore, we have
\[
\chi_i \ = \ \omega_f^{E_i}
\ = \ \omega_f^{(p^f-1)/(p-1)}
\ = \ \chi^{\mathrm{cyc}}_p,
\]
as desired.
\end{proof}

\begin{lemma}[Cyclotomic rigidity]\label{lem:ff-cyclotomic}
We continue to use the notation of Lemma \ref{lem:raynaud}.
Suppose that, for any $i\in\{1,\ldots,r\}$, $\chi_i=\chi^{\mathrm{cyc}}_p$.
Then $\mathcal{G}$ is isomorphic to $\mu_p^r$ as an $R$-group scheme.
\end{lemma}

\begin{proof}
Write
\[
\mathcal E \ \stackrel{\mathrm{def}}{=} \ \mathcal G^\vee
\]
for the Cartier dual of $\mathcal G$, and
\[
W \ \stackrel{\mathrm{def}}{=} \ \mathcal E_K(\overline K).
\]
By Cartier duality, we have a natural isomorphism of $\mathbb F_p[G_K]$-modules
\[
W \ \stackrel{\sim}{\rightarrow} \ \operatorname{Hom}_{\mathbb F_p} (V,\mu_p(\overline K)) \ \ (\stackrel{\sim}{\rightarrow} \ \operatorname{Hom}_{\mathbb F_p} (V,\mathbb{F}_p)\otimes_{\mathbb{F}_p}\mu_p(\overline K)).
\]
By extending scalars to $\overline{\mathbb F}_p$, it follows that the characters associated to the Jordan--H\"older constituents of $W\otimes_{\mathbb{F}_p}\overline{\mathbb{F}}_p$ are given by
\[
\chi_p^{\mathrm{cyc}}\chi_i^{-1} \qquad (1 \ \le \ i \ \le \ r),
\]
which are trivial by our assumption.

Now we consider the connected--\'etale exact
sequence of $\mathcal E$ [cf. \cite{TateFFGS}, \S 3.7, (I)]:
\[
0 \ \longrightarrow \ \mathcal E^\circ \ \longrightarrow \ \mathcal E \ \longrightarrow \ \mathcal E^{\mathrm{et}} \ \longrightarrow \ 0,
\]
where $\mathcal E^\circ$ is a finite flat closed subgroup scheme of $\mathcal E$
with connected special fiber, and $\mathcal E^{\mathrm{et}}$ is a finite \'etale group scheme over $R$.
I claim that $\mathcal E^\circ$ is trivial.
Indeed, suppose that $\mathcal E^\circ$ is nontrivial, hence that [since $\mathcal E^\circ$ is finite flat over $R$ and its generic fiber is finite \'etale over $K$] $\mathcal E^\circ_K(\overline K) \neq \{0\}$.
Then since $\mathcal E^\circ \subseteq \mathcal E$ is a closed subgroup scheme,
$\mathcal E^\circ_K(\overline K)$ may naturally be regarded as a nontrivial $\mathbb{F}_p[G_K]$-submodule of $W$.
In particular, the characters associated to the Jordan--H\"older constituents of $\mathcal E^\circ_K(\overline K)\otimes_{\mathbb{F}_p}\overline{\mathbb{F}}_p$ are trivial.
On the other hand, since $(\mathcal E^\circ)_k$ is connected, it follows from Lemma \ref{lem:raynaud}, (ii), that these characters are nontrivial, a contradiction.
Thus, $\mathcal E^\circ$ is trivial.

In light of this claim, $\mathcal E$ is finite \'etale over $R$.
Since $R$ is strictly henselian, $\mathcal E$ is a constant group scheme.
Here, we note that $\mathcal E$ is annihilated by $p$ and has rank $p^r$ [cf. Lemma \ref{lem:raynaud}, (i)].
Therefore,
\[
\mathcal E \ \simeq \ (\mathbb Z/p\mathbb Z)_R^r.
\]
Taking Cartier duals, we conclude that
\[
\mathcal G \ \simeq \ \mu_p^r,
\]
as desired.
\end{proof}

\begin{proposition}[Bounded-order multiplicative rigidity]\label{prop:ff-bounded}
We continue to use the notation of Lemma \ref{lem:raynaud}.
Suppose that
\begin{itemize}
\item $\mathcal G_k$ is connected, and
\item each $\chi_i$ has order at most $p(p-1)$.
\end{itemize}
Then $\mathcal G$ is isomorphic to $\mu_p^r$ as an $R$-group scheme.
\end{proposition}

\begin{proof}
The assertion follows from Lemmas \ref{lem:digit} and \ref{lem:ff-cyclotomic}.
\end{proof}

\section{Displacements and their finite flat closure}\label{sec:displacement}
In the present section, we study the Galois displacements of a point on a curve under the assumption that these displacements are torsion.
By using an auxiliary finite flat group scheme, we show that the group generated by these displacements is annihilated by $p$.
Then we construct the finite flat closure of this group and show that its special fiber is connected.

In the present section, suppose that we are in the situation at the beginning of \S 1. Let
\begin{itemize}
\item $g\ge2$ be an integer,
\item $X$ a smooth, proper, geometrically connected curve over $K$ of genus $g$, and
\item $x \in X(\overline{K})$.
\end{itemize}
Write
\begin{itemize}
\item $S \stackrel{\mathrm{def}}{=} G_Kx$ ($\subseteq X(\overline{K})$),
\item $m$ for the cardinality of $S$,
\item $J$ for the Jacobian variety of $X$, and
\item $D \stackrel{\mathrm{def}}{=} \big\langle[y-x] \mid y \in S\big\rangle_{\mathbb{Z}}$ ($\subseteq J(\Kbar)$).
\end{itemize}
For each proper group scheme $\mathcal Z$ over $R$, write
\[
\operatorname{red}_{\mathcal Z}:\mathcal Z_K(\Kbar) \ \rightarrow \ \mathcal Z_k(k)
\]
for the natural reduction homomorphism.
We note that if, moreover, $\mathcal Z$ is finite flat over $R$, then $\operatorname{red}_{\mathcal Z}$ is surjective [cf. our assumption that $k$ is algebraically closed and $R$ is henselian; \cite{BLR}, Chapter 2, \S 2.3, Proposition 4].

Suppose that 
\begin{itemize}
\item $J$ has good reduction over $R$, i.e., there exists an abelian scheme $\mathcal{J}$ over $R$ equipped with an isomorphism $\mathcal{J}_K \simeq J$, and
\item for each $\sigma\in G_K$, $[\sigma x-x]\in J(\Kbar)$ is torsion.
\end{itemize}

\begin{lemma}[Finite flat realization of division points]
\label{lem:finite-flat-kummer}
Let $n\ge 2$ be an integer, and
\[
b \ \in \ J(K) \ \ (\stackrel{\sim}{\rightarrow} \ \cJ(R)).
\]
Then there exists a finite flat group scheme $\cH_b$ over $R$ satisfying the following conditions:
\begin{itemize}
\item $\cH_b$ is commutative,
\item $\cH_b$ is annihilated by $n$,
\item there exists an exact sequence of $R$-group schemes
\[
0 \ \longrightarrow \ \cJ[n]
\ \stackrel{\iota}{\longrightarrow} \ \cH_b
\ \stackrel{q}{\longrightarrow} \ (\mathbb Z/n\mathbb Z)_R
\ \longrightarrow \ 0,
\]
and
\item for each $a\in J(\Kbar)$ with $na=b$ [which implies that --- since $b \in J(K)$ --- $\sigma a-a \in \cJ[n](\Kbar)$ for each $\sigma\in G_K$], there exists an element
\[
u \ \in \ \cH_b(\Kbar)
\]
such that
\[
\sigma u-u \ = \ \iota(\sigma a-a)
\]
for each $\sigma\in G_K$.
\end{itemize}
\end{lemma}

\begin{proof}
For each $i\in\{0,\ldots,n-1\}$, write
\[
\cH_i \ \stackrel{\mathrm{def}}{=} \ [n]^{-1}(ib) \ \ (\subseteq \ \cJ)
\]
for the scheme-theoretic inverse image of the section $ib \in \cJ(R)$ via the multiplication-by-$n$ morphism
\[
[n]:\cJ \ \rightarrow \ \cJ.
\]
Since $[n]$ is finite flat, each $\cH_i$ is a finite flat scheme over $R$.
Now we consider the finite flat scheme
\[
\cH_b \ \stackrel{\mathrm{def}}{=} \ \coprod_{i=0}^{n-1}\cH_i
\]
over $R$.
We equip $\cH_b$ with a commutative group scheme structure over $R$ as follows.
Let $T$ be an $R$-scheme, and
\[
z\in\cH_i(T),\qquad w\in\cH_j(T),
\]
where $i,j\in\{0,\ldots,n-1\}$.
Write $(z,i)$ and $(w,j)$ for the corresponding $T$-valued points of $\cH_b$.
Define
\[
(z,i)+(w,j)
\ \stackrel{\mathrm{def}}{=}\
\begin{cases}
(z+w,i+j), & i+j<n,\\
(z+w-b_T,i+j-n), & i+j\ge n,
\end{cases}
\]
where $b_T\in\cJ(T)$ denotes the image of $b\in\cJ(R)$.
[We note that since
\[
n(z+w-b_T) \ = \ (i+j-n)b_T,
\]
these formulas are well-defined.]
In particular, these formulas define an $R$-morphism
\[
\cH_b\times_R\cH_b \ \rightarrow \ \cH_b.
\]
One verifies directly that this operation defines a commutative group scheme structure on $\cH_b$ with identity $(0,0)$ and inverses given by
\[
-(z,i)
\ =\
\begin{cases}
(-z,0), & i=0,\\
(b_T-z,n-i), & i>0.
\end{cases}
\]
Thus, $\cH_b$ is a finite flat commutative group scheme over $R$.
Moreover, it follows from the equalities
\[
n(z,i) \ = \ (nz-ib_T,0) \ = \ (0,0)
\]
that $\cH_b$ is annihilated by $n$.

Write
\[
q:\cH_b \ = \ \coprod_{i=0}^{n-1}\cH_i \ \rightarrow \ \coprod_{i=0}^{n-1}\operatorname{Spec}(R) \ = \ (\mathbb Z/n\mathbb Z)_R
\]
for the natural projection, and
\[
\iota:\cJ[n] \ = \ \cH_0 \ \hookrightarrow \ \cH_b
\]
for the natural inclusion.
By the definition of the group law, $q$ is a homomorphism of $R$-group schemes with kernel
\[
\ker(q) \ = \ \cH_0 \ = \ \cJ[n].
\]
Here, we note that since $[n]$ is finite faithfully flat, each $\cH_i$ is a finite faithfully flat scheme over $R$.
Thus, it follows that $q$ is finite faithfully flat, hence the sequence
\[
0 \ \longrightarrow \ \cJ[n]
\ \stackrel{\iota}{\longrightarrow} \ \cH_b
\ \stackrel{q}{\longrightarrow} \ (\mathbb Z/n\mathbb Z)_R
\ \longrightarrow \ 0
\]
is exact.

Finally, let $a\in J(\Kbar)$ be an element such that $na=b$.
If we write
\[
u \ \stackrel{\mathrm{def}}{=} \ (a,1) \ \in \ \cH_{b}(\Kbar),
\]
then for each $\sigma\in G_K$, we have
\[
\sigma u-u \ = \ (\sigma a,1)+(b-a,n-1) \ = \ (\sigma a-a,0) \ = \ \iota(\sigma a-a).
\]
\end{proof}

\begin{lemma}[Exponent-$p$ rigidity of Galois displacements]
\label{lem:kh-descent}
It holds that
\[
pD \ = \{0\}.
\]
In particular, $D$ is a $\mathbb{F}_p[G_K]$-submodule of $J[p](\Kbar)$.
Moreover, if $(X,x)$ is exceptional [cf. \cite{Tamagawa}, Definition preceding
Proposition 1.4], then 
\[
D \ = \ \{0\},
\]
hence that $x\in X(K)$ [cf. our assumption that $g\ge2$].
\end{lemma}

\begin{proof}
We begin by claiming that 
\[
X(K) \ \ne \ \emptyset.
\]
Indeed, we note that since $J$ has good reduction over $R$, $X$ has stable reduction over $R$ [cf. \cite{DeligneMumford}, Theorem 2.4].
Write $\mathcal X$ for the stable model of $X$ over $R$, and $\mathcal U$ for its smooth locus over $R$.
Then it follows from the fact that $\mathcal X_k$ is a stable curve that
\[
\mathcal U_k(k) \ \ne \ \emptyset
\]
[cf. our assumption that $k$ is algebraically closed].
Thus, by \cite{BLR}, Chapter 2, \S 2.3,
Proposition 5, we obtain that
\[
\mathcal{U}(R) \ \ne \ \emptyset,
\]
hence that $X(K)\ne\emptyset$. This completes the proof of the claim.

Fix an element $x_0\in X(K)$.
Write
\[
a \ \stackrel{\mathrm{def}}{=} \ [x-x_0] \ \in \ J(\Kbar).
\]
Here, observe that since the natural action of $G_K$ on $\cJ_k(k)$ is trivial [cf. our assumption that $k$ is algebraically closed], for each $\sigma\in G_K$, we have
\[
\operatorname{red}_{\cJ}(\sigma a-a)
 \ = \ \sigma\operatorname{red}_{\cJ}(a)-\operatorname{red}_{\cJ}(a)
 \ = \ 0.
\]
Moreover, since $\operatorname{red}_{\cJ}$ is injective on $J(\Kbar)[p^{'\infty}]$ [cf. \cite{SerreTate1968}, \S 1, Lemma 2], it follows from our assumption that $\sigma a-a=[\sigma x-x]$ is torsion that 
\[
\sigma a-a \ \in \ J(\Kbar)[p^\infty].
\]
Thus, there exists an integer $s\ge1$ such that if we write $n\stackrel{\mathrm{def}}{=}p^s$, then, for each $\sigma \in G_K$, we have
\[
n(\sigma a-a) \ = \ 0
\]
[cf. the fact that $S$ is a finite set].
In particular, it holds that
\[
b \ \stackrel{\mathrm{def}}{=} \ na \ \in J(K).
\]

First, we consider the case where $(X,x)$ is exceptional.
Then it follows from the proof of \cite{Tamagawa}, Proposition 3.1, (i), that, for each $\sigma\in G_K$, we have
\[
2\sigma x \ \sim \ 2x.
\]
In particular, we obtain that
\[
2(\sigma a-a) \ = \ 0.
\]
Thus, it follows from the fact that $\sigma a-a$ is $p$-primary [and $p$ is odd] that
\[
\sigma a-a \ = \ 0,
\]
hence that $D=\{0\}$.

Next, we consider the case where $(X,x)$ is not exceptional.
Observe that, by Lemma \ref{lem:finite-flat-kummer}, there exists a finite flat group scheme $\cH_b$ over $R$ satisfying the following conditions:
\begin{itemize}
\item $\cH_b$ is commutative,
\item $\cH_b$ is annihilated by $n$,
\item there exists an injective $G_K$-equivariant homomorphism $\iota:\cJ[n](\Kbar)\rightarrow\cH_b(\Kbar)$, and
\item there exists an element $u\in\cH_b(\Kbar)$ such that
\[
\sigma u-u \ = \ \iota(\sigma a-a)
\]
for each $\sigma\in G_K$.
\end{itemize}
Write
\[
M \ \subseteq \ \cH_b(\Kbar)
\]
for the $\mathbb{Z}[G_K]$-submodule generated by $u$.
We note that this $G_K$-module $M$ is
\begin{itemize}
\item finite [cf. the fact that $\cH_b$ is finite flat over $R$],
\item $p$-primary [i.e., $M = M_p$ --- cf. the fact that $\cH_b$ is annihilated by $n=p^s$], and
\item of type $G$, hence of type $S$ [cf. \cite{Hoshi2017}, Definition 2.3, (i), and Remark 2.3.1, (i), (iii)].
\end{itemize}

Now I claim that
\begin{itemize}
\item $u$ is a weak $G_K$-invariant, i.e., for any $\gamma,\delta\in G_K$, the following implication holds:
\[
(1-\gamma)^2(\delta u) \ = \ 0 \quad \Longrightarrow \quad (1-\gamma)(\delta u) \ = \ 0.
\]
\end{itemize}
Indeed, let $\gamma,\delta\in G_K$ be elements such that
\[
(1-\gamma)^2(\delta u) \ = \ 0.
\]
Then since
\[
\sigma u-u \ = \ \iota(\sigma a-a)
\]
for each $\sigma\in G_K$, we have
\[
0 \ = \ (1-\gamma)^2(\delta u) \ = \ \iota(\delta a-2\gamma\delta a+\gamma^2\delta a).
\]
[Here, we note that since $na=b\in J(K)$,
$\delta a-2\gamma\delta a+\gamma^2\delta a
\in\cJ[n](\Kbar)$.]
Thus, it follows from the injectivity of $\iota$ that
\[
\delta a-2\gamma\delta a+\gamma^2\delta a \ = \ 0,
\]
hence that
\[
\delta x+\gamma^2\delta x \ \sim \ 2\gamma\delta x
\]
[cf. the fact that $a=[x-x_0]$ and $x_0\in X(K)$].
In particular, since $(X,\gamma\delta x)$ is not exceptional [cf. our assumption that $(X,x)$ is not exceptional], we obtain an equality of divisors
\[
\delta x+\gamma^2\delta x \ = \ 2\gamma\delta x.
\]
Therefore, we conclude that
\[
\delta x \ = \ \gamma\delta x,
\]
hence that
\[
(1-\gamma)(\delta u)
\ = \ \iota(\delta a-\gamma\delta a)
\ = \ 0.
\]
This completes the proof of the claim.

In light of this claim [together with the above three itemized properties of $M$], it follows from \cite{Hoshi2017}, Lemma 2.7, that the action of
$G_K$ on $pM_p=pM$ is trivial.
In particular, for each $\sigma\in G_K$, we have
\[
p(\sigma u-u) \ = \ 0.
\]
Therefore, since
\[
p(\sigma u-u)
\ = \ \iota(p(\sigma a-a))
\]
and $\iota$ is injective, we conclude that
\[
p(\sigma a-a) \ = \ 0,
\]
hence that $pD = \{0\}$.
\end{proof}

\begin{lemma}[Finite flat displacement group with connected special fiber]\label{lem:closure}
Write $\cD$ for the schematic closure in $\cJ[p]$ of the closed subgroup scheme of $\cJ[p]_K$ corresponding to $D$ [cf. Lemma \ref{lem:kh-descent}]. 
Then $\cD$ is a finite flat closed subgroup scheme of $\cJ[p]$ over $R$ such that 
\begin{itemize}
\item $\cD$ is commutative,
\item $\cD$ is annihilated by $p$, and
\item $\cD_k$ is connected.
\end{itemize}
\end{lemma}

\begin{proof}
Since $\cJ$ is an abelian scheme over $R$, $\cJ[p]$ is a finite flat closed subgroup scheme of $\cJ$ over $R$ such that 
\begin{itemize}
\item $\cJ[p]$ is commutative, and
\item $\cJ[p]$ is annihilated by $p$.
\end{itemize}
In particular, Lemma \ref{lem:ff-closure} implies that $\cD$ is a finite flat closed subgroup scheme of $\cJ[p]$ over $R$ such that 
\begin{itemize}
\item $\cD$ is commutative, and
\item $\cD$ is annihilated by $p$.
\end{itemize}

Let us verify the connectedness of $\cD_k$.
We observe that it follows from the discussion in the second paragraph of the proof of Lemma \ref{lem:kh-descent} that $\operatorname{red}_{\cJ}(D)=\{0\}$.
On the other hand, since the natural reduction homomorphism
\[
\operatorname{red}_{\cD}:D \ = \ \cD_K(\Kbar) \ \rightarrow \ \cD_k(k)
\]
is surjective, and the natural homomorphism $\cD_k(k) \rightarrow \cJ_k(k)$ is injective, we have 
\[
\cD_k(k) \ = \ \{0\}.
\]
Therefore, we conclude that $\cD_k$ is a one-point scheme, hence connected.
\end{proof}

\section{Orbit augmentation and a dimension inequality}\label{sec:k0}
In the present section, we establish a bound on the dimension of the space generated by Galois displacements in terms of the size of the Galois orbit.

In the present section, suppose that we are in the situation at the beginning of \S 2. Write
\begin{itemize}
\item $P_S \stackrel{\mathrm{def}}{=} \Fp[S]=\oplus_{y \in S}\Fp e_y$, where $G_K$ acts on $P_S$ by $\sigma e_y=e_{\sigma y}$, 
\item $I_S \subseteq P_S$ for the kernel of the homomorphism $P_S \rightarrow \Fp$ given by sending
\[
\sum_{y \in S} n_y e_y \ \mapsto \ \sum_{y \in S} n_y,
\]
where $n_y \in \Fp$, and
\item $T \subseteq I_S$ for the kernel of the homomorphism $\rho: I_S \rightarrow D$ given by sending 
\[
\sum_{y \in S} n_y e_y \ \mapsto \ \sum_{y\in S}n_y[y-x].
\]
\end{itemize}
Here, we note that, for any $\sigma\in G_K$ and $u=\sum_{y\in S}n_y e_y\in I_S$, we have
\[
\sigma\rho(u) \ = \ \sum_{y\in S}n_y[\sigma y-\sigma x] \ = \ \sum_{y\in S}n_y[\sigma y-x] - \Bigl(\sum_{y\in S}n_y\Bigr)[\sigma x-x] \ = \ \rho(\sigma u),
\]
where the final equality follows from the fact that $\sum_{y\in S}n_y=0$.
Thus, $\rho$ is $G_K$-equivariant.
Also, note that, for each $y\in S$, we have
\[
\rho(e_y-e_x)=[y-x].
\]
Since $D$ is generated by the elements $[y-x]$ for $y\in S$,
we conclude that $\rho$ is surjective.
In summary, we obtain an exact sequence of $\Fp[G_K]$-modules
\[
0 \ \longrightarrow \ T \ \longrightarrow \ I_S \ \stackrel{\rho}{\longrightarrow} \ D \ \longrightarrow \ 0.
\]
  
\begin{proposition}[Orbit augmentation inequality]\label{prop:k0}
We have an inequality
\[
\dim_{\Fp}D \ \le \ \dim_{\Fp}T.
\]
\end{proposition}

\begin{proof}
We begin by observing the following elementary fact:
Let $G$ be a group, $F$ a field, and $B$ a finite $G$-set.
Then we have a natural isomorphism of $F[G]$-modules
\[
F[B] \ \stackrel{\sim}{\rightarrow} \ F[B]^\vee,
\]
where we write $(-)^\vee$ for the dual of the $F[G]$-module $(-)$.
In light of this fact, we have natural isomorphisms of $\Fp[G_K]$-modules
\[
P_S \ \stackrel{\sim}{\rightarrow} \ P^\vee_S, \qquad \Fp \ \stackrel{\sim}{\rightarrow} \ \Fp^\vee.
\]
In particular, the augmentation exact sequence
\[
0 \ \longrightarrow \ I_S \ \longrightarrow \ P_S \ \longrightarrow \ \Fp \ \longrightarrow \ 0
\]
and its dual imply that
\[
[I_S] \ =\ [P_S]-[\Fp] \ = \ [P_S^\vee]-[\Fp^\vee] \ = \ [I_S^\vee]
\]
in the Grothendieck group of finite-dimensional $\Fp[G_K]$-modules.

For each finite-dimensional $\Fp[G_K]$-module $V$, write
\[
\overline{V} \ \stackrel{\mathrm{def}}{=} \ V\otimes_{\Fp}\overline{\mathbb{F}}_p.
\]
For each character $\lambda:G_K \rightarrow \overline{\mathbb{F}}^{\times}_p$, write
\[
m_{\overline{V}}(\lambda)
\]
for the Jordan--Hölder multiplicity of $\lambda$ in $\overline{V}$.
In particular, we have the following equalities:
\[
m_{\overline{I}_S}(\lambda) \ = \ m_{\overline{I}_S}(\lambda^{-1}), \qquad
m_{\overline{I}_S}(\lambda) \ = \ m_{\overline{T}}(\lambda)+m_{\overline{D}}(\lambda),
\]
where [since scalar extension is exact and commutes with duals]
the first equality follows from the equality $[I_S] = [I_S^\vee]$, and the second equality follows from the exact sequence
\[
0 \ \longrightarrow \ T \ \longrightarrow \ I_S \ \stackrel{\rho}{\longrightarrow} \ D
\ \longrightarrow \ 0.
\]

Write $\Sigma$ for the set of distinct characters occurring as Jordan--Hölder constituents of $\overline{D}$.
Since $\mathcal D_k$ is connected [cf. Lemma \ref{lem:closure}], it follows from
Lemma \ref{lem:raynaud}, (i), (iii), that each Jordan--Hölder constituent of $\overline{D}$ is one-dimensional and that
\[
\Sigma \cap \Sigma^{-1} \ =\ \varnothing.
\]
Thus, for each $\chi\in\Sigma$, we obtain that $m_{\overline{D}}(\chi^{-1})=0$, hence that
\[
m_{\overline{T}}(\chi^{-1}) \ = \ m_{\overline{I}_S}(\chi^{-1}) \ = \ m_{\overline{I}_S}(\chi) \ = \ m_{\overline{T}}(\chi)+m_{\overline{D}}(\chi) \ \ge \ m_{\overline{D}}(\chi).
\]
Therefore, since the characters $\chi^{-1}$ for $\chi\in\Sigma$ are pairwise distinct, summing over $\Sigma$ yields
\[
\dim_{\Fp}D \ = \ \dim_{\overline{\mathbb{F}}_p}\overline{D} \ = \ \sum_{\chi\in\Sigma}m_{\overline{D}}(\chi) \ \le \ \sum_{\chi\in\Sigma}m_{\overline{T}}(\chi^{-1}) \ \le \ \dim_{\overline{\mathbb{F}}_p}\overline{T} \ = \ \dim_{\Fp}T.
\]
\end{proof}

\section{Kummer theory and the exclusion of pencils}\label{sec:pencil}
In the present section, we combine the character rigidity established in \S 1, the dimension inequality established in \S 3, and Kummer theory to obtain a lower bound for
the dimension of the space of global sections associated to $px$.

In the present section, suppose that we are in the situation at the beginning of \S 3. Write
\begin{itemize}
\item $X_{\Kbar} \stackrel{\mathrm{def}}{=} X\times_{K}\Kbar$, 
\item $\mathcal{O}(E)$ for the line bundle associated to a divisor $E$ on $X_{\Kbar}$,
\item $h^0(E) \stackrel{\mathrm{def}}{=} \dim_{\Kbar}H^0(X_{\Kbar}, \mathcal{O}(E))$, and
\item $H$ for the image of $G_K$ in the permutation group $\operatorname{Sym}(S)$ of $S$. 
\end{itemize}

\begin{lemma}[Kummer restriction bound]\label{lem:kummer}
Suppose that $m>1$ and $h^0(px)=2$. 
Then we have an inequality
\[
\dim_{\Fp}T \ \le \ 1.
\]
\end{lemma}

\begin{proof}
Since $m>1$, fix $y_0\in S\setminus\{x\}$.
Then by Lemma \ref{lem:kh-descent}, there exists a nonconstant rational function $f\in K(X_{\Kbar})^\times$ such that
\[
\divisor(f) \ = \ p(y_0-x),
\]
where we write $K((-))$ for the function field of the integral scheme $(-)$.
Since $1,f\in H^0(X_{\Kbar},\mathcal O(px))$ are linearly independent and $h^0(px)=2$, we have
\[
H^0(X_{\Kbar},\mathcal O(px)) \ = \ \Kbar\cdot1\oplus \Kbar\cdot f.
\]
We note that, in a similar vein, for each $z\in S\setminus\{x\}$, there exists a nonconstant rational function $g_z \in K(X_{\Kbar})^\times$ such that
\[
\divisor(g_z) \ = \ p(z-x).
\]
In particular, since $g_z\in H^0(X_{\Kbar},\mathcal O(px))$, there exist $a_z\in \Kbar$ and $b_z\in \Kbar^\times$ such that
\[
g_z \ = \ a_z+b_z f,
\]
which implies that $f(z)=-a_z/b_z$.
Thus, if we write [by slight abuse of notation] 
\[
f:X_{\Kbar} \ \to \ \bP^1_{\Kbar}
\]
for the morphism determined by the rational function $f$, then we have
\begin{itemize}
\item $\deg(f)=p$,
\item $f(x)=\infty$,
\item $f(y_0)=0$, and
\item $f(z)=-a_z/b_z$ for each $z\in S\setminus\{x\}$.
\end{itemize}
Here, we observe that $f$ is injective on $S$.
[Indeed, since $f(x)=\infty$ and $f(z)\in\Kbar$ for each $z\in S\setminus\{x\}$, it suffices to consider $z,w\in S\setminus\{x\}$.
Suppose that $f(z)=f(w)$.
Then it follows that $g_z/b_z=g_w/b_w$, hence that $\divisor(g_z)=\divisor(g_w)$.
Thus, we obtain that $z=w$.]

Write
\begin{itemize}
\item $B \stackrel{\mathrm{def}}{=} f(S)\subseteq\bP^1(\Kbar)$,
\item $\operatorname{Div}^0(\bP^1_{\Kbar})$ for the group of degree-zero divisors on $\bP^1_{\Kbar}$, and
\item $\operatorname{Div}_B^0(\bP^1_{\Kbar})$ for the group of degree-zero divisors on $\bP^1_{\Kbar}$ supported on $B$.
\end{itemize}
Then the bijection $f|_S:S\stackrel{\sim}{\to}B$ induces an isomorphism of $\Fp$-vector spaces
\[
I_S \ \stackrel{\sim}{\rightarrow} \ \operatorname{Div}_B^0(\bP^1_{\Kbar})/
p\operatorname{Div}_B^0(\bP^1_{\Kbar});
\quad
\sum_{z\in S}(\widetilde n_z\ (\mathrm{mod}\ p))e_z \ \mapsto\ \sum_{z\in S}\widetilde n_z f(z) \ (\mathrm{mod}\ p\operatorname{Div}_B^0(\bP^1_{\Kbar})),
\]
where $\widetilde n_z\in\mathbb Z$ are chosen so that $\sum_{z\in S}\widetilde n_z=0$.
On the other hand, it follows from the equality $\mathrm{Pic}^0(\bP^1_{\Kbar})=\{0\}$ that we have an exact sequence of abelian groups 
\[
0 \ \longrightarrow \ \Kbar^{\times} \ \longrightarrow \ K(\bP^1_{\Kbar})^{\times} \ \stackrel{\mathrm{div}}{\longrightarrow} \ \operatorname{Div}^0(\bP^1_{\Kbar})
\ \longrightarrow \ 0.
\]
Thus, we obtain an isomorphism of $\Fp$-vector spaces
\[
K(\bP^1_{\Kbar})^{\times}/K(\bP^1_{\Kbar})^{\times p} \ \stackrel{\sim}{\rightarrow} \ 
\operatorname{Div}^0(\bP^1_{\Kbar})/p\operatorname{Div}^0(\bP^1_{\Kbar})
\]
[cf. the fact that $\Kbar$ is algebraically closed].
In summary, we have a composite of injective $\Fp$-linear maps
\[
\kappa: I_S \ \stackrel{\sim}{\rightarrow} \ \operatorname{Div}_B^0(\bP^1_{\Kbar})/
p\operatorname{Div}_B^0(\bP^1_{\Kbar}) \ \hookrightarrow \ \operatorname{Div}^0(\bP^1_{\Kbar})/p\operatorname{Div}^0(\bP^1_{\Kbar}) \ \stackrel{\sim}{\rightarrow} \ K(\bP^1_{\Kbar})^{\times}/K(\bP^1_{\Kbar})^{\times p},
\]
where the second arrow is induced by the natural injection $\operatorname{Div}_B^0(\bP^1_{\Kbar}) \hookrightarrow \operatorname{Div}^0(\bP^1_{\Kbar})$.

Now I claim that
\[
(T \ \stackrel{\sim}{\rightarrow}) \ \ \kappa(T) \ \subseteq \ \ker\left(
K(\bP^1_{\Kbar})^\times/K(\bP^1_{\Kbar})^{\times p} \ \stackrel{K(f)}{\rightarrow} \ K(X_{\Kbar})^\times/K(X_{\Kbar})^{\times p}
\right).
\]
Indeed, we first observe that
\begin{itemize}
\item $f^*(\infty)=px$ [cf. the equality $\divisor(f) = p(y_0-x)$],
\item $f^*(0)=py_0$ [cf. the equality $\divisor(f) = p(y_0-x)$], and
\item $f^*(f(z))=pz$ for each $z\in S\setminus\{x\}$ [cf. the equalities $\divisor(f-f(z)) = \divisor(g_z/b_z) = \divisor(g_z) = p(z-x)$ and $f^*(\infty)=px$].
\end{itemize}
Let $t=\sum_{z\in S}n_z e_z\in T$.
Choose integer lifts $\widetilde n_z$ such that $\sum_{z\in S}\widetilde{n}_z=0$.
Write
\[
E_S=\sum_{z\in S}\widetilde n_z z, \qquad
E_B=\sum_{z\in S}\widetilde n_z f(z).
\]
Then it follows from the equality $\mathrm{Pic}^0(\bP^1_{\Kbar})=\{0\}$ that
there exists a rational function $q\in K(\bP^1_{\Kbar})^\times$ such that 
\[
\divisor(q) \ = \ E_B.
\]
By the fiber equalities above, we obtain that
\[
\divisor(q\circ f) \ = \ f^*\divisor(q) \ = \ f^*E_B \ = \ f^*\Bigl(\sum_{z\in S}\widetilde n_z f(z)\Bigr) \ = \ \sum_{z\in S}\widetilde n_z pz \ = \ pE_S.
\]
On the other hand, since $t\in T=\ker\rho$ and $p[z-x]=0$ for each $z\in S$, we have 
\[
[E_S] \ = \ \Bigl[\sum_{z\in S}\widetilde n_z z\Bigr] \ = \ \sum_{z\in S}\widetilde n_z [z-x] \ = \ \rho(t) \ = \ 0 \ \in \ J(\Kbar).
\]
Thus, there exists a rational function $u\in K(X_{\Kbar})^\times$ such that 
\[
\divisor(u) \ = \ E_S.
\]
Therefore, we conclude that 
\[
(q\circ f)/u^p \ \in \ \Kbar^{\times} \ = \ \Kbar^{\times p} \ \subseteq \ K(X_{\Kbar})^{\times p},
\]
hence that $\kappa(t)$ is trivial in $K(X_{\Kbar})^\times/K(X_{\Kbar})^{\times p}$.
This completes the proof of the claim.

Write
\begin{itemize}
\item $L \stackrel{\mathrm{def}}{=} K(\bP^1_{\Kbar})$,
\item $M \stackrel{\mathrm{def}}{=} K(X_{\Kbar})$, and
\item $\Gamma \subseteq L^{\times}$ for the inverse image of $\kappa(T)$ via the natural surjection $L^\times\to L^\times/L^{\times p}$.
\end{itemize}
Since $L$ contains a primitive $p$-th root of unity ($\in \Kbar$), the above claim implies that
\[
L(\Gamma^{1/p}) \ \subseteq \ M.
\]
Thus, it follows from Kummer theory [cf., e.g., \cite{Lang-Algebra}, Chapter VI, Section 8, Theorem 8.1] that
\[
[L(\Gamma^{1/p}):L] \ = \ [\Gamma:L^{\times p}] \ = \ p^{\dim_{\Fp}T}.
\]
On the other hand, since
\[
[L(\Gamma^{1/p}):L] \ \mid \ [M:L] \ = \ \deg(f) \ = \ p, 
\]
we conclude that
\[
p^{\dim_{\Fp}T} \ \mid \ p,
\]
hence that $\dim_{\Fp}T\le1$.
\end{proof}

\begin{proposition}[Lower bound for $h^0(px)$]
\label{prop:no-pencil}
Suppose that $m>1$.
Then we have an inequality
\[
h^0(px) \ \ge \ 3.
\]
\end{proposition}

\begin{proof}
We note that it follows from the discussion at the beginning of the proof of Lemma \ref{lem:kummer} that $h^0(px)\ge2$.
Thus, to verify Proposition \ref{prop:no-pencil}, it suffices to show that $h^0(px)\neq2$.

Suppose that $h^0(px)=2$.
Then it follows from Proposition \ref{prop:k0} and Lemma \ref{lem:kummer} that
\[
\dim_{\Fp}D \ \le \ \dim_{\Fp}T \ \le \ 1.
\]
On the other hand, since $m>1$, it holds that $D\ne\{0\}$ [cf. our assumption that $g\ge2$], hence that
\[
\dim_{\Fp}D \ = \ \dim_{\Fp}T \ = \ 1.
\]
Thus, we obtain that
\[
m-1 \ = \ \dim_{\Fp}I_S \ = \ \dim_{\Fp}T+\dim_{\Fp}D \ = \ 2,
\]
hence that $m=3$.

Since 
\begin{itemize}
\item the representation $G_K \rightarrow \operatorname{Aut}_{\Fp}(I_S)$ factors as the composite of the natural homomorphism $G_K \rightarrow \operatorname{Sym}(S) = \operatorname{Sym}(G_Kx)$ and the natural homomorphism $\operatorname{Sym}(S) \rightarrow \operatorname{Aut}_{\Fp}(I_S)$, and
\item $\rho:I_S \rightarrow D$ is a $G_K$-equivariant surjection,
\end{itemize}
it follows that the representation $\chi_D: G_K \rightarrow \operatorname{Aut}_{\Fp}(D)$ factors as the composite of the natural surjection $G_K \twoheadrightarrow H$ and the natural homomorphism $\delta: H \rightarrow \operatorname{Aut}_{\Fp}(D)$.
Moreover, since $\dim_{\Fp}D = 1$, $\delta$ factors as the composite of the natural surjection $H \twoheadrightarrow H^{\mathrm{ab}}$ and the natural homomorphism 
\[
\overline{\delta}: H^{\mathrm{ab}} \ \rightarrow \ \operatorname{Aut}_{\Fp}(D).
\]

Here, we note that since $H$ acts faithfully and transitively on $S$ and $m=3$, it holds that
\[
H \ \simeq \ \mathbb{Z}/3\mathbb{Z} \qquad \text{or} \qquad H \ \simeq \ \operatorname{Sym}(S),
\]
hence that the order of $H^{\mathrm{ab}}$ is at most $3$.
In particular, the character $\chi_D: G_K \rightarrow \operatorname{Aut}_{\Fp}(D) \simeq \Fp^{\times}$ has order at most $3$.
On the other hand, since $\cD_k$ is connected [cf. Lemma \ref{lem:closure}],
it follows from Lemma \ref{lem:digit} that 
\[
\chi_D \ = \ \chi^{\mathrm{cyc}}_p.
\]
This contradicts the fact that $\chi^{\mathrm{cyc}}_p$ has order $p-1\ge4$.
Therefore, we conclude that $h^0(px)\neq2$.
\end{proof}

\section{Small-orbit exclusion}\label{sec:small}
In the present section, we establish a lower bound for the size of a
nontrivial Galois orbit by combining Guralnick's theorem on cyclic
quotients of transitive groups, multiplicative rigidity for finite
flat group schemes, and the geometry of degree-two divisors on curves.

In the present section, suppose that we are in the situation at the beginning of \S 4.
Write
\begin{itemize}
\item $d \stackrel{\mathrm{def}}{=} \dim_{\Fp}D$, and
\item $\chi_D: G_K \rightarrow \operatorname{Aut}_{\Fp}(D)$
for the representation associated to the natural action of $G_K$ on $D$.
\end{itemize}

\begin{lemma}[Multiplicative rigidity for small orbits]\label{lem:mu}
Suppose that $m\le p(p-1)$.
Then $\mathcal D$ is isomorphic to $\mu_p^{d}$ as an $R$-group scheme.
\end{lemma}

\begin{proof}
Since 
\begin{itemize}
\item the representation $G_K \rightarrow \operatorname{Aut}_{\Fp}(I_S)$ factors as the composite of the natural homomorphism $G_K \rightarrow \operatorname{Sym}(S) = \operatorname{Sym}(G_Kx)$ and the natural homomorphism $\operatorname{Sym}(S) \rightarrow \operatorname{Aut}_{\Fp}(I_S)$, and
\item $\rho:I_S \rightarrow D$ is a $G_K$-equivariant surjection,
\end{itemize}
it follows that $\chi_D$ factors through $H$.
Write
\[
\chi_1,\ldots,\chi_d: H \ \rightarrow \ \overline{\mathbb{F}}^{\times}_p
\]
for the characters associated to the Jordan--Hölder constituents of $D\otimes_{\Fp}\overline{\mathbb{F}}_p$ [cf. Lemma \ref{lem:raynaud}, (i), and Lemma \ref{lem:closure}].
For each $i\in\{1,\ldots,d\}$, write
\[
N_i \ \stackrel{\mathrm{def}}{=} \ \ker(\chi_i).
\]
Then since $H$ acts transitively on $S$ and $H/N_i$ is cyclic of order $\operatorname{ord}(\chi_i)$ [cf. the fact that each finite subgroup of the
multiplicative group of a field is cyclic], it follows from \cite{Guralnick2000}, Theorem 1, that
\[
\operatorname{ord}(\chi_i) \ \le \ m \ \le \ p(p-1).
\]

In light of the inequalities $\operatorname{ord}(\chi_i) \le p(p-1)$, Lemma \ref{lem:mu} follows from the fact that $\mathcal D_k$ is connected [cf. Lemma \ref{lem:closure}] and Proposition \ref{prop:ff-bounded}.
\end{proof}

\begin{lemma}[Cyclotomic displacement rigidity]\label{lem:cocycle}
Suppose that $G_K$ acts on $D$ via the mod-$p$
cyclotomic character.
Then $x\in X(K)$.
\end{lemma}

\begin{proof}
Suppose that $x \in X(\Kbar)\setminus X(K)$, and hence that 
\begin{itemize}
\item $(X,x)$ is not exceptional [cf. Lemma \ref{lem:kh-descent}], and
\item $D\neq\{0\}$ [cf. our assumption that $g\ge2$].
\end{itemize}
Let $\sigma \in G_K$ be an element such that
\[
\chi_p^{\mathrm{cyc}}(\sigma) \ = \ -2
\]
[cf. the surjectivity of $\chi_p^{\mathrm{cyc}}$].
For each $z\in S$, write
\[
d_z \ \stackrel{\mathrm{def}}{=} \ [\sigma z-z] \ = \ [\sigma z-x]-[z-x] \ \in \ D.
\]
Then since $G_K$ acts on $D$ via $\chi_p^{\mathrm{cyc}}$, it holds that
\[
[\sigma^2z-\sigma z]
\ = \ \sigma d_z
\ = \ -2d_z
\ = \ -2[\sigma z-z],
\]
hence that
\[
\sigma^2z+\sigma z \ \sim \ 2z.
\]
In particular, since $(X,z)$ is not exceptional [cf. our assumption that $(X,x)$ is not exceptional], we obtain an equality of divisors
\[
\sigma^2z+\sigma z \ = \ 2z.
\]
This implies that $\sigma z=z$.
Thus, $\sigma$ acts trivially on $S$, hence on $D$.
This contradicts the fact that $\sigma$ acts on the nonzero $\Fp$-vector space $D$ as multiplication by $-2\neq1$ [cf. our assumption that $p\ge5$].
Therefore, we conclude that $x\in X(K)$.
\end{proof}

\begin{proposition}[Common small-orbit exclusion]\label{prop:small}
Suppose that $m>1$.
Then $m> p(p-1)$.
\end{proposition}

\begin{proof}
Suppose that $m\le p(p-1)$.
Then it follows from Lemma \ref{lem:mu} that we have an isomorphism of
$\mathbb F_p[G_K]$-modules
\[
D \ \simeq \ \mu_p(\overline K)^d.
\]
In light of Lemma \ref{lem:cocycle}, this implies that
\[
x \ \in \ X(K),
\]
hence that $m=1$, a contradiction.
\end{proof}

\section{A geometric bound for prime-degree linear systems}\label{sec:geometry}
In the present section, independently of the arithmetic
setting of the preceding sections, we establish a purely
geometric bound for divisors on smooth proper curves over algebraically closed fields of characteristic zero.
Our main ingredients are the Pl\"ucker formula, Clifford's
theorem, Castelnuovo's bound, and the Riemann--Roch theorem.

\begin{proposition}[Prime-degree geometric bound]\label{prop:prime-degree}
Let
\begin{itemize}
\item $p\ge5$ be a prime number,
\item $g\ge2$ an integer,
\item $F$ an algebraically closed field of characteristic zero,
\item $C$ a smooth, proper, connected curve over $F$ of genus $g$,
\item $E$ a divisor of degree $p$ on $C$, and
\item $B\subseteq C(F)$ a subset such that
      $py\sim E$ for each $y\in B$.
\end{itemize}
[Note that $B$ is automatically finite. Indeed, write $J$ for the Jacobian variety of $C$. If $B\ne\emptyset$ and $y_0\in B$, then the image of $B$ via the Abel--Jacobi map
\[
C(F) \ \rightarrow \ J(F); \quad y \ \mapsto \ [y-y_0]
\]
is a subset of $J[p](F)$. Since this map is injective, we have $|B|\le p^{2g}$.]
Suppose that
\[
\dim_F H^0(C,\mathcal O(E)) \ \ge \ 3,
\]
where we write $\mathcal{O}(E)$ for the line bundle associated to $E$. 
Then we have an inequality
\[
|B| \ \le \ (p-1)^2.
\]
\end{proposition}

\begin{proof}
If $|B|\le1$, then the assertion is immediate.
Thus, we may assume that $|B|>1$, and hence that the complete linear
system $|E|$ is base-point-free. 
[Indeed, if we take two distinct points $y,z\in B$, then $py,pz\in|E|$, and these divisors have disjoint supports.]
Write
\begin{itemize}
\item $r \stackrel{\mathrm{def}}{=} \dim_F H^0(C,\mathcal O(E))-1 \ge 2$,
\item $\phi:C \rightarrow \mathbb P_F^r$ for the morphism associated to $|E|$,
\item $\Gamma \stackrel{\mathrm{def}}{=} \phi(C)
      \subseteq \mathbb P_F^r$,
      endowed with the reduced induced closed
      subscheme structure, and
\item $\overline{\phi}:C \rightarrow \Gamma$
      for the [surjective] morphism induced by $\phi$.
\end{itemize}
Here, observe that since
\begin{itemize}
\item $\Gamma$ is nondegenerate,
\item irreducible curves in $\mathbb P_F^r$ of degree one are lines, and
\item lines in $\mathbb P_F^r$ are degenerate [cf. the fact that $r\ge2$],
\end{itemize}
we have $\deg(\Gamma)>1$. 
Therefore, we conclude from the degree formula
\[
p \ = \ \deg(E) \ = \ \deg(\overline{\phi})\deg(\Gamma)
\]
that
\[
\deg(\overline{\phi}) \ = \ 1, \qquad \deg(\Gamma) \ = \ p
\]
[cf. our assumption that $p$ is a prime number].

For each $z\in C(F)$, write
\[
0 \ \le \ h_0(z) \ < \ \cdots \ < \ h_r(z) \ \le \ p
\]
for the integers $n\ge0$ satisfying
\[
h^0(E-nz) \ > \ h^0(E-(n+1)z),
\]
where we write $h^0((-))$ for the dimension over $F$ of the space of global sections of the line bundle associated to the divisor $(-)$.
[One verifies easily that these integers coincide with the Hermite invariants of $H^0(C,\mathcal O(E))$ at $z$ defined in \cite{Laksov1984}, \S 4, by considering a basis with strictly increasing orders of vanishing.]
In particular, for each $y \in B$, since $py\sim E$, we obtain that
\[
h^0(E-py) \ = \ 1,\qquad h^0(E-(p+1)y) \ = \ 0,
\]
hence that $h_r(y)=p$.
For each $z\in C(F)$, write
\[
d_r(z) \ \stackrel{\mathrm{def}}{=} \ \sum_{i=0}^r(h_i(z)-i).
\]
Then, for each $y \in B$, we have
\[
d_r(y) \ = \ \sum_{i=0}^r(h_i(y)-i) \ = \ (h_r(y)-r) + \sum_{i=0}^{r-1}(h_i(y)-i) \ \ge \ p-r
\]
[cf. the fact that $h_i(y)\ge i$].
Thus, it follows from the Pl\"ucker formula [cf. \cite{Laksov1984}, Theorem 9 and
Theorem 15, (i) and (ii); also the definition preceding Theorem 9] that
\begin{equation}\label{eq:Plucker}
|B|(p-r) \ \le \ \sum_{y\in B}d_r(y) \ \le \ \sum_{z\in C(F)}d_r(z) \ = \ (r+1)(p+r(g-1)).
\end{equation}

First, suppose that $p<2g$.
Note that since $|B|>1$, there exists $y \in C(F)$ such that $E \sim py$.
Then it follows from Clifford's theorem [cf. \cite{ACGH1985}, Chapter III, \S 1,
Clifford's Theorem, (i)] that
\[
r \ = \ h^0(E)-1 \ = \ h^0(py)-1 \ \le \ \frac{p}{2}
\]
[cf. the fact that $py$ is effective], hence that
\[
2r+1 \ \le \ p
\]
[cf. our assumption that $p$ is odd].
Thus, we have
\begin{equation}\label{eq:Clifford}
\frac{r+1}{p-r} \ \le \ 1.
\end{equation}
On the other hand, since $\overline{\phi}$ is a birational morphism onto the nondegenerate curve $\Gamma$ of degree $p$ in $\mathbb{P}^r_F$, it follows from Castelnuovo's bound
[cf. \cite{ACGH1985}, Chapter III, \S 2, Castelnuovo's Bound] that
\[
g \ \le \ \frac{m(m-1)}{2}(r-1)+m\varepsilon,
\]
where $(m, \varepsilon)$ is the pair of nonnegative integers such that
\[
p-1 \ = \ m(r-1)+\varepsilon, \qquad 0 \ \le \ \varepsilon\ < \ r-1.
\]
Then since
\[
\frac{(p-1)(p-2)}{2(r-1)}-\Bigl\{\frac{m(m-1)}{2}(r-1)+m\varepsilon\Bigr\} \ = \ \frac{m(r-2)}{2}
+\frac{\varepsilon(\varepsilon-1)}{2(r-1)} \ \ge \ 0
\]
[cf. the fact that $r\ge2$], we obtain that
\begin{equation}\label{eq:Castelnuovo}
g \ \le \ \frac{(p-1)(p-2)}{2(r-1)}
\end{equation}
Therefore, we conclude that
\begin{align*}
|B| \
&\le \ \frac{r+1}{p-r}(p+r(g-1))\\
&\le \ p+r(g-1)\\
&\le \ p+\frac{r(p-1)(p-2)}{2(r-1)}-r\\
&\le \ p+(p-1)(p-2)-2 \ = \ p^2-2p \ < \ (p-1)^2,
\end{align*}
where 
\begin{itemize}
\item the first inequality follows from \eqref{eq:Plucker},
\item the second inequality follows from \eqref{eq:Clifford},
\item the third inequality follows from \eqref{eq:Castelnuovo}, and
\item the fourth inequality follows from $r\ge2$.
\end{itemize}

Next, suppose that $p>2g$.
Write $K_C$ for the canonical divisor of $C$.
Then we observe that
\[
\deg(K_C-E) \ = \ 2g-2-p \ < \ 0,
\]
hence that 
\[
h^0(K_C-E) \ = \ 0.
\]
Thus, it follows from the Riemann--Roch theorem that
\[
r+1 \ = \ h^0(E) \ = \ p-g+1,
\]
hence that
\[
(r+1)(p+r(g-1)) \ = \ g(p-g+1)^2.
\]
Therefore, we conclude from \eqref{eq:Plucker} that
\[
|B|g \ \le \ g(p-g+1)^2,
\]
hence that
\[
|B| \ \le \ (p-g+1)^2 \ \le \ (p-1)^2
\]
[cf. our assumption that $g\ge2$].
\end{proof}

\section{Proof of the main theorem}\label{sec:main}

In the present section, we prove our main theorem by combining the Kummer-theoretic exclusion of pencils in \S 4, the small-orbit exclusion in \S 5, and the geometric bound established in \S 6.

In the present section, suppose that we are in the situation at the beginning of \S 5.

\begin{theorem}[The local Coleman assertion]\label{thm:main}
It holds that
\[
m \ = \ 1.
\]
Therefore, we conclude that
\[
x \ \in \ X(K).
\]
\end{theorem}

\begin{proof}
Suppose that $m>1$.
Then it follows from Proposition \ref{prop:small} that
\[
m \ > \ p(p-1).
\]
Now we consider the divisor $px$ on $X_{\Kbar}$.
Then by Lemma \ref{lem:kh-descent}, each element $y \in S$ satisfies $py \sim px$.
Moreover, by Proposition \ref{prop:no-pencil}, we have
\[
h^0(px) \ \ge \ 3.
\]
Thus, it follows from Proposition \ref{prop:prime-degree} that
\[
m \ = \ |S| \ \le \ (p-1)^2 \ < \ p(p-1),
\]
contradicting $m>p(p-1)$. Therefore, we conclude that $m=1$.
\end{proof}

\medskip
\noindent\textbf{Acknowledgments.}
This research was supported by 
JSPS KAKENHI Grant Number JP24K16898,
the ZEN Mathematics Center at ZEN University, and
the Research Institute for Mathematical Sciences, an International Joint Usage/Research Center located in Kyoto University. 
This research was also supported by the CNRS France-Japan ``Arithmetic and Homotopic Galois Theory'' International Research Network.

\medskip
\noindent\textbf{Human--AI collaboration.}
This work was developed through extensive mathematical discussions
between the author and OpenAI's ChatGPT.
AI assistance included mathematical exploration, the development
and examination of arguments, reference checking, and drafting
and revising the manuscript.
The author takes full responsibility for the mathematical
content of the paper.

\end{document}